\documentclass[11pt,amsfonts]{article}
\usepackage{graphicx}
\usepackage{latexsym}
\usepackage{amssymb}
\usepackage{amsmath}
\usepackage{layout}
\newtheorem{prop}{Proposition}
\newtheorem{lemma}{Lemma}

\newtheorem{corollary}{Corollary}
\newtheorem{theorem}{Theorem}
\newtheorem{remark}{Remark}

\def\real{{\mathord{{\rm I\kern-2.8pt R}}}}        
\def\inte{{\mathord{{\rm I\kern-2.8pt N}}}}

\def\sZZ{{\rm Z\kern-2.8ptem{}Z}}

\def\z{{\mathchoice
  {\sZZ}
  {\sZZ}
  {\rm Z\kern-0.30em{}Z}
  {\rm Z\kern-0.25em{}Z} }}
\def\sQQ{{\kern 0.27em \vrule height1.45ex width0.03em depth0em
          \kern-0.30em \rm Q}}
\def\qu{{\mathchoice
    {\sQQ}
    {\sQQ}
  {\kern 0.225em \vrule height1.05ex width0.025em depth0em \kern-0.25em \rm Q}
  {\kern 0.180em \vrule height0.78ex width0.020em depth0em \kern-0.20em \rm Q}
        }}
\def\sCC{{\kern 0.27em \vrule height1.45ex width0.03em depth0em
          \kern-0.30em \rm C}}
\def\complex{{\mathchoice
    {\sCC}
    {\sCC}
  {\kern 0.225em \vrule height1.05ex width0.025em depth0em \kern-0.25em \rm C}
  {\kern 0.180em \vrule height0.78ex width0.020em depth0em \kern-0.20em \rm C}
        }}

\newcommand{\ba}{\begin{array}}
\newcommand{\ea}{\end{array}}
\newcommand{\be}{\begin{equation}}
\newcommand{\ee}{\end{equation}}
\newcommand{\bea}{\begin{eqnarray}}
\newcommand{\eea}{\end{eqnarray}}
\newcommand{\beaa}{\begin{eqnarray*}}
\newcommand{\eeaa}{\end{eqnarray*}}

\def\z{\zeta}

\font\tenmath=msbm10 \font\sevenmath=msbm7 \font\fivemath=msbm5
\newfam\mathfam \textfont\mathfam=\tenmath
\scriptfont\mathfam=\sevenmath \scriptscriptfont\mathfam=\fivemath

\def \={{\buildrel {\rm (law)} \over =}}

\def\qed{ \hfill \vrule width.25cm height.25cm depth0cm\smallskip}

\newcommand{\basa}{\begin{assumption}}
\newcommand{\easa}{\end{assumption}}

\newcommand{\bas}{\begin{assum}}
\newcommand{\eas}{\end{assum}}

\newcommand{\ignore}[1]{}
\begin{document}

\renewcommand{\thefootnote}{\fnsymbol{footnote}}

\renewcommand{\thefootnote}{\fnsymbol{footnote}}

\date{ August 19, 2009 }
\title{On the structure of Gaussian random variables}
\author{Ciprian A. Tudor\\
  SAMOS/MATISSE,
Centre d'Economie de La Sorbonne,\\ Universit\'e de
Panth\'eon-Sorbonne Paris 1,\\90, rue de Tolbiac, 75634 Paris Cedex
13, France.\\tudor@univ-paris1.fr\vspace*{0.1in}} \maketitle

\begin{abstract}
We study when a given  Gaussian random variable on a given  probability space $\left( \Omega , {\cal{F}}, P\right) $
 is equal almost surely to $\beta_{1}$ where $\beta $ is a Brownian motion defined on the same
  (or possibly extended) probability space. As a consequence of this result, we prove that the distribution of a random variable
  in a  finite sum of Wiener chaoses  (satisfying in addition a certain property) cannot be normal. This result also allows to understand better a characterization of the Gaussian variables obtained via Malliavin calculus.
\end{abstract}

\vskip0.5cm

{\bf  2000 AMS Classification Numbers:} 60G15,  60H05,  60H07.

 \vskip0.3cm


\section{Introduction}
We study when a Gaussian random variable defined on some  probability space can be expressed almost surely as a Wiener integral
 with respect to a Brownian motion defined on the same space. The starting point of this work are some recent results related to
 the distance between the law of  an arbitrary random variable $X$ and the Gaussian law. This distance can be defined in various ways
 (the Kolmogorov distance, the total variations distance or others) and it can be expressed in terms of the Malliavin derivative $DX$
 of  the random variable $X$ when this derivative exists.  These results lead to a characterization of Gaussian random variables through
  Malliavin calculus. Let us briefly recall the context.
Suppose that $\left( \Omega, {\cal{F}}, P\right)$ is  a probability space and let $(W_{t})_{t\in[0,1]}$ be a ${\cal{F}}_{t}$  Brownian motion
 on this space, where ${\cal{F}}_{t}$ is its natural filtration. Equivalent conditions for the standard normality of a centered
  random variable $X$ with variance $1$   are the following:
  $\mathbf{E}\left( 1-\langle DX, D(-L)^{-1} \rangle| X\right) =0$ or $\mathbf{E}\left( f'_{z}(X) (1-  \langle DX, D(-L)^{-1} \rangle\right)=0$
   for every $z$ where $D$ denotes the Malliavin derivative, $L$ is the Ornstein-Uhlenbeck operator,
    $\langle \cdot , \cdot \rangle $ denotes the scalar product in $L^{2}([0,1])$  and the deterministic function
    $f_{z}$ is the solution of the Stein's equation (see e.g. \cite{NoPe1}).
    This characterization is of course interesting and it can be useful in some cases.
    It is also easy to understand it for random variables that are Wiener integrals with respect to $W$.
    Indeed, assume that  $X=W(h)$ where $h$ is a deterministic function in $L^{2}([0,1]) $  with $\Vert h\Vert _{L^{2}([0,1])}=1$.
    In this case $DX= h= D(-L)^{-1} X$ and then $\langle DX, D(-L)^{-1} \rangle=1$ and the above equivalent conditions for the
     normality of $X$  can be easily verified.  In some other cases, it is difficult, even impossible, to compute the quantity
     $\mathbf{E}\left(\langle DX, D(-L)^{-1} \rangle| X\right) $ or $\mathbf{E}\left( f'_{z}(X) (1-  \langle DX, D(-L)^{-1} \rangle\right)$.
      Let us consider for example the case of the random variable $Y=\int_{0}^{1} sign (W_{s} ) dW_{s}$.
      This is not a Wiener integral with respect to $W$. But it is  well-known that it is standard  Gaussian because
       the process $\beta_{t}= \int_{0}^{t} sign(W_{s}) dW_{s}$ is a Brownian motion as follows from the L\'evy's characterization theorem.
       The chaos expansion of this random variable is known  and it is recalled in Section \ref{2}.
       In fact  $Y$ can be expressed as an infinite sum of multiple Wiener-It\^o stochastic integrals and it is impossible to check
       if the equivalent conditions for its normality are satisfied (it is even not differentiable in the Malliavin calculus  sense).
        The phenomenon that happens here is that $Y$ can be expressed as the value at time 1 of the Brownian motion $\beta $ which is
        actually the Dambis-Dubbins-Schwarz (DDS in short) Brownian motion associated to the martingale $M^{Y}= (M^{Y}_{t})_{t\in [0,1] }$,
        $M^{Y}_{t}=\mathbf{E}\left( Y | {\cal{F}}_{t}\right) $ (recall that ${\cal{F}}_{t}$ is the natural filtration  of $W$  and
        $\beta $ is defined on the same  space $\Omega$ (or possibly on a extension of $\Omega$) and is a ${\cal{G}}_{s} $ -Brownian motion
        with respect to the filtration ${\cal{G}}_{s}= {\cal{F}}_{T(s)}$ where $T(s)= \inf ( t\in [0,1]; \langle M^{Y}\rangle _{t} \geq s)$).
This leads   to  the following question: is any standard normal
random variable $X$ representable as the value at time 1 of the
Brownian motion associated, via the Dambis-Dubbins-Schwarz theorem,
to the martingale $M^{X}$,  where for every $t$
\begin{equation}\label{mx}
M^{X} _{t} =\mathbf{E}(X| {\cal{F}}_{t}) ?
\end{equation}
By combining the techniques of Malliavin calculus and classical
tools of the probability theory, we found the following answer: if
the bracket of the ${\cal{F}}_{t}$ martingale $M^{X}$ is bounded
a.s. by 1 then this property is true, that is, $X$ can be
represented as its DDS Brownian motion at time 1. The property also
holds when the bracket $\langle M^{X} \rangle _{1} $ is bounded by
an arbitrary constant and $\langle M^{X} \rangle _{1} $ and $\beta
_{\langle M^{X} \rangle _{1}} $ are independent. If the bracket of
$M^{X}$ is not bounded by 1, then this property is not true. An
example when it fails is obtained by considering the standard normal
random variable $W(h_{1})sign(W(h_{2})) $  where $h_{1}, h_{2}$ are
two orthonormal elements of $L^{2}([0,1])$. Nevertheless, we will
prove that we can construct a bigger probability space $\Omega _{0}$
that includes $\Omega$ and a Brownian motion on $\Omega _{0}$ such
that $X$ is equal almost surely with this Brownian motion at time 1.
The construction is done by the means of the Karhunen-Lo\`eve
theorem.
 Some consequences of this result are discussed here; we believe that these consequences could be various.
  We prove that the standard normal random variables such that the bracket of its associated DDS martingale
  is bounded by 1 cannot live in a finite sum of Wiener chaoses: they can be or in the first chaos, or in an infinite sum of chaoses.
   We also make a connection with some results obtained recently via Stein's method and Malliavin calculus.

 We structured our paper as follows. Section 2 starts with a short description of the elements of the
 Malliavin calculus and it also contains   our main result on the structure of Gaussian random variables.
 In Section 3 we discusses some consequences of our characterization. In particular we prove that the random variables
 whose associated DDS martingale has  bracket bouned by 1 cannot belong to a finite sum of Wiener chaoses and we relate our work with recent results on standard normal random variables obtained via Malliavin calculus.

\section{On the structure of Gaussian random variable}
\label{2}
Let us consider a probability space $\left( \Omega, {\cal{F}} , P\right)$ and assume that $(W_{t})_{t\in [0,1]} $ is a
Brownian motion on this space with respect to its natural  filtration $\left( {\cal{F}}_{t} \right) _{t\in [0,1]}$.  Let $I_{n}$ denote
 the multiple Wiener-It\^o  integral of order $n$  with respect to $W$.
 The elements of the stochastic calculus for multiple integrals and of Malliavin calculus can be found in \cite{Mal} or \cite{N}.
 We will just introduce  very briefly some notation.
    We recall that any square integrable random variable which is
measurable with respect to the $\sigma$-algebra generated by $W$ can
be expanded into an orthogonal sum of multiple stochastic integrals
\begin{equation}
\label{sum1} F=\sum_{n\geq0}I_{n}(f_{n})
\end{equation}
where $f_{n}\in L^{2}([0,1]^{n})$ are (uniquely determined)
symmetric functions and $I_{0}(f_{0})=\mathbf{E}\left[  F\right]  $.

The isometry of multiple integrals can be written as: for $m,n$ positive integers and $f\in L^{2}([0,1]^{n})$, $g\in L^{2}([0,1] ^{m})$

\begin{eqnarray}
\mathbf{E}\left(I_{n}(f) I_{m}(g) \right) &=& n! \langle f,g\rangle _{L^{2}([0,1])^{\otimes n}}\quad \mbox{if } m=n,\nonumber \\
\mathbf{E}\left(I_{n}(f) I_{m}(g) \right) &= & 0\quad \mbox{if } m\not=n.\label{iso}
\end{eqnarray}
It also holds that
\begin{equation*}
I_{n}(f) = I_{n}\big( \tilde{f}\big)
\end{equation*}
where $\tilde{f} $ denotes the symmetrization of $f$ defined by $\tilde{f}%
(x_{1}, \ldots , x_{x}) =\frac{1}{n!} \sum_{\sigma \in {\cal S}_{n}}
f(x_{\sigma (1) }, \ldots , x_{\sigma (n) } ) $.
We will need the general formula for calculating products of Wiener chaos
integrals of any orders $m,n$ for any symmetric integrands $f\in
L^{2}([0,1]^{\otimes m})$ and $g\in L^{2}([0,1]^{\otimes n})$; it is%
\begin{equation}
I_{m}(f)I_{n}(g)=\sum_{l=0}^{m\wedge n}l!C_{m}^{l}C_{n}^{l}%
I_{m+m-2l}(f\otimes_{l}g) \label{product}%
\end{equation}
where the contraction $f\otimes_{l}g$ ($0\leq l\leq m\wedge n$) is
defined by
\begin{eqnarray}
&&  (f\otimes_{\ell} g) ( s_{1}, \ldots, s_{n-\ell}, t_{1}, \ldots, t_{m-\ell
})\nonumber\\
&&  =\int_{[0,T] ^{m+n-2\ell} } f( s_{1}, \ldots, s_{n-\ell}, u_{1},
\ldots,u_{\ell})g(t_{1}, \ldots, t_{m-\ell},u_{1}, \ldots,u_{\ell})
du_{1}\ldots du_{\ell} . \label{contra}%
\end{eqnarray}
Note that the contraction $(f\otimes_{\ell} g) $ is an element of $L^{2}([0,1]^{m+n-2\ell})$ but it is not necessary symmetric. We will by $(f\tilde{\otimes }_{\ell} g)$ its symmetrization.

We denote by $\mathbb{D}^{1,2}$ the domain of the Malliavin
derivative with respect to $W$ which takes values in $L^{2}([0,1]
\times \Omega)$. We just recall that $D$ acts on functionals of the
form $f(X)$, with $X\in \mathbb{D}^{1,2}$ and $f$ differentiable,
in the following way: $D_{\alpha} f(X)= f'(X)D_{\alpha }X$ for every
$\alpha \in (0,1]$  and on multiple integrals $I_{n}(f)$ with $f\in
L^{2}([0,1] ^{n})$ as $D_{\alpha}I_{n}(f)= nI_{n-1} f(\cdot ,
\alpha)$.

 The Malliavin derivative $D$ admits a dual operator which is the divergence integral $\delta(u)\in L^{2}(\Omega)$ if $u\in Dom(\delta)$ and we have  the duality relationship
\begin{equation}\label{dua}
\mathbf{E}(F\delta (u))= \mathbf{E}\langle DF, u\rangle, \hskip0.5cm F\in \mathbb{D}^{1,2}, u\in Dom(\delta).
\end{equation}
For adapted integrands, the divergence integral coincides with
the classical It\^o integral.

\vskip0.2cm

Let us fix the probability space $(\Omega ,{\cal{F}}, P)$ and let us
assume that the Wiener process $(W_{t})_{t\in [0,1]} $ lives on this
space. Let $X$ be a centered square integrable random variable on
$\Omega$. Assume that $X$ is measurable with respect to the
sigma-algebra ${\cal{F}}_{1}$. After Proposition \ref{p1} the random
variable  $X$ will be assumed to have standard normal law.

The following result is an immediate consequence of the Dambis-Dubbins-Schwarz theorem (DDS theorem for short, see \cite{KS}, Section 3.4, or \cite{RY}, Chapter V).
\begin{prop}\label{p1}
Let $X$ be a random variable in $L^{1}(\Omega)$. Then there exists a Brownian motion $(\beta _{s})_{s\geq 0}$  (possibly defined on an extension of the probability space) with respect to a filtration $({\cal{G}}_{s})_{s\geq 0}$ such that
\begin{equation*}
X= \beta _{\langle M^{X} \rangle _{1}}
\end{equation*}
where $M^{X}=(M^{X}_{t} )_{t\in [0,1]}$ is the martingale given by (\ref{mx}). Moreover the random time $T=\langle M^{X}\rangle _{1} $ is a stopping time for the filtration ${\cal{G}}_{s}$ and it satisfies $T> 0$ a.s. and $\mathbf{E}T=\mathbf{E}X^{2}.$
\end{prop}
{\bf Proof: } Let $T(s)=\inf  \left( t\geq 0,  \langle M^{X} \rangle _{t} \geq s\right)$. By applying Dambis-Dubbins-Schwarz theorem
\begin{equation*}
\beta_{s}:= M_{T(s)}
\end{equation*}
is a standard Brownian motion with respect to the filtration ${\cal{G}}_{s}:={\cal{F}}_{T(s)}$ and for every $t\in [0,1]$ we have $M^{X}_{t}= \beta _{\langle M^{X}\rangle _{t}}$ a.s. $P$. Taking $t=1$ we get
\begin{equation*}
X= \beta _{\langle M^{X} \rangle _{1}} \hskip0.5cm \mbox{ a.s.}.
\end{equation*}

The fact that $T$ is a $({\cal{G}}_{s})_{s\geq 0}$ stopping time is well known. It is true because $(\langle M^{X}\rangle _{1} \leq s)= (T(s)\geq 1) \in {\cal{F}}_{T(s)}={\cal{G}}_{s}$. Also clearly $T>0$ a.s and $\mathbf{E}T= \mathbf{E} X^{2}$.\qed

\vskip0.3cm

In the sequel we will  call the Brownian $\beta$ obtained via the DDS theorem as the DDS Brownian associated to $X$.

Recall the Ocone-Clark formula: if $X$ is a random variable in $\mathbb{D}^{1,2}$ then
\begin{equation}
\label{oc}
X=\mathbf{E}X+ \int_{0}^{1} \mathbf{E}\left( D_{\alpha} X| {\cal{F}}_{\alpha }\right) dW_{\alpha}.
\end{equation}

\begin{remark}
If the random variable $X$ has zero mean and it belongs to the space $\mathbb{D}^{1,2}$ then by the Ocone-Clark formula (\ref{oc}) we have $M^{X}_{t}= \int_{0}^{t} \mathbf{E}\left( D_{\alpha }X| {\cal{F}}_{\alpha } \right) dW_{\alpha }$
and consequently
\begin{equation*}
X= \beta _{\int_{0}^{1} \left( \mathbf{E}\left( D_{\alpha }X| {\cal{F}}_{\alpha } \right)\right) ^{2} d\alpha}
\end{equation*}
where $\beta $ is the DDS Brownian motion associated to $X$.
\end{remark}

Assume from now on   that $X\sim N(0,1)$. As we have seen, $X$ can be written as the value at a random time of a Brownian motion $\beta$ (which is fact the Dambis-Dubbins-Schwarz Brownian associated to the martingale $M^{X}$). Note that $\beta$ has the time interval $\mathbb{R}_{+}$ even if $W$  is indexed over $[0,1]$. So, if we know that $\beta_{T}$ has a standard normal law, what can we say about the random time $T$? It is equal to 1 almost surely? This is for example  the case of the variable $X= \int_{0}^{1} sign(W_{s}) dW_{s}$ because here, for every $t\in [0,1]$,  $M^{X} _{t}=\int_{0}^{t} sign(W_{s})dW_{s} $ and $\langle M^{X}\rangle _{t}= \int_{0}^{t} (sign(B_{s})^{2}ds =t$. An other situation when this is true is related to Bessel processes.   Let  $(B^{(1)}, \ldots B^{(d)}) $ be a $d$-dimensional Brownian motion and consider   the random variable
\begin{equation}\label{bes}
X=\int_{0}^{1} \frac{B^{(1)}_{s}}{\sqrt{ (B^{(1)}_{s} )^{2}+ \ldots + (B^{(d)}_{s} )^{2}}}dB^{(1)}_{s}+ \ldots +\int_{0}^{1} \frac{B^{(d)}_{s}}{\sqrt{ (B^{(1)}_{s} )^{2}+ \ldots + (B^{(d)}_{s} )^{2}}}dB^{(d)}_{s}
\end{equation}
It also satisfies $T:=\langle M^{X} \rangle _{t}=t$ for every $t\in
[0,1]$ and in particular  $\langle M^{X} \rangle _{1}=1$ a.s.. We
will see below that the fact that any $N(0,1)$ random variable is
equal a.s. to
 $\beta _{1}$ (its associated DDS Brownian evaluated at time 1) is true only for random variables for
 which the bracket of their  associated DDS martingale is almost surely bounded and $T$ and $\beta _{T}$ are independent
 or if $T$ is bounded  almost surely by $1$.

\vskip0.2cm

We will assume the following condition on the stopping time $T$.

\begin{equation}
\label{H}
\mbox{ There exist a constant $M>0$ such that $T\leq M$ a.s. }
\end{equation}

The problem we address  in this   section is then the following: let $(\beta _{t})_{t\geq 0}$ be a ${\cal{G}}_{t}$- Brownian motion
 and let $T$ be a almost surely positive stopping time for its filtration such that $\mathbf{E}(T)=1$ and $T$ satisfies (\ref{H}).
 We will show when $T=1$ a.s.

\vskip0.2cm

Let us start with the following result.

\begin{theorem}
Assume (\ref{H}) and assume that $T$ is independent by $\beta _{T}$.
Then it holds that $\mathbf{E}T^{2}=1$.
\end{theorem}
{ \bf Proof: } Let us apply It\^o's formula to the ${\cal{G}}_{t}$ martingale $\beta_{T\wedge t}$. Letting $t\to \infty $  (recall that $T$ is a.s. bounded) we get
\begin{equation*}
\mathbf{E}\beta_{T} ^{4} = 6\mathbf{E}\int_{0}^{T} \beta_{s}^{2}ds. \end{equation*}

Since $\beta_{T}$ has $N(0,1)$ law, we have that $\mathbf{E}\beta_{T}^{4}=3$. Consequently
\begin{equation*}
\mathbf{E}\int_{0}^{T} \beta_{s}^{2}ds=\frac{1}{2}.
\end{equation*}
Now, by the independence of $T$ and $\beta _{T}$, we get
$\mathbf{E}(T\beta_{T}^{2})=\mathbf{E}T \mathbf{E}\beta_{T}^{2}=1$.
Applying again It\^o formula to $\beta_{T\wedge t}$ with $f(t,x)=
tx^{2}$ we get
\begin{equation*}
\mathbf{E}T\beta_{T}^{2} =\mathbf{E}\int_{0}^{T} \beta_{s}^{2}ds +
\mathbf{E}\int_{0}^{T} sds.
\end{equation*}
Therefore $\mathbf{E}\int_{0}^{T} sds=\frac{1}{2}$ and then $\mathbf{E}T^{2}=1$. \qed

\begin{theorem}\label{tt}
Let $(\beta_{t})_{t\geq 0}$ be  a ${\cal{G}}_{t}$ Wiener process and
let $T$ be a ${\cal{G}}_{t}$ bounded stopping time with
$\mathbf{E}T=1$. Assume that $T$ and $\beta _{t}$ are independent.
Suppose $\beta_{T}$ has a $N(0,1)$ law. Then $T=1$ a.s.
\end{theorem}
{\bf Proof: } It is a consequence of the above proposition, since
$\mathbf{E}(T-1) ^{2}= \mathbf{E}T^{2}-2\mathbf{E}(T)+1=0$. \qed

\vskip0.3cm

\begin{prop}\label{bounded1} Assume that $(\ref{H})$ is satisfied with $M \leq 1$. Then $T=1$ almost surely.
\end{prop}
{\bf Proof: } By It\^o's formula,
\begin{equation*}
 \mathbf{E} \beta _{T}^{4} = 6\mathbf{E}\int_{0}^{T} \beta
_{s}^{2}ds=6\mathbf{E}\int_{0}^{1} \beta _{s}^{2}ds+
\mathbf{E}\int_{\mathbb{R}_{+}} \beta _{s}^{2}1_{[T,1]}(s)ds .
\end{equation*}
Since $6\mathbf{E}\int_{0}^{1} \beta _{s}^{2}ds=3$ and $ \mathbf{E}
\beta _{T}^{4} =3$ it follows that $\mathbf{E}\int_{\mathbb{R}_{+}}
\beta _{s}^{2}1_{[T,1]}(s)ds =0$ and this implies that $ \beta
_{s}^{2}(\omega)1_{[T(\omega) ,1]}(s)=0$ for almost all $s$ and
$\omega$. Clearly $T=1$ almost surely. \qed


\vskip0.3cm

Next, we will try  to understand if this property is always true without the assumption that the bracket of the martingale $M^{X}$ is finite almost surely. To this end, we will consider the following example. Let $(W_{t})_{t\in [0,1]} $ a standard  Wiener process with respect to its natural filtration ${\cal{F}}_{t}$. Consider $h_{1}, h_{2}$ two functions in $L^{2}([0,1])$ such that $\langle h_{1}, h_{2} \rangle _{L^{2}([0,1])}=0$ and $\Vert h_{1}\Vert _{L^{2}([0,1])}= \Vert h_{2}\Vert _{L^{2}([0,1])}=1$. For example we can choose
$$h_{1}(x)=\sqrt{2} 1_{[0,\frac{1}{2}]}(x) \mbox{  and  }h_{2}(x)= \sqrt{2} 1_{[\frac{1}{2}, 1]}(x)$$ (so, in addition, $h_{1}$ and $h_{2}$ have disjoint support). Define the random variable
\begin{equation}\label{x}
X=  W(h_{1}) sign W(h_{2}).
\end{equation}
It is well-known that $X$ is standard normal. Note in particular that $X^{2} =W(h_{1}) ^{2}$. We will see that it cannot be written as the value at time 1 of its associated DDS martingale.  To this end we will use the chaos expansion of $X$ into multiple Wiener-It\^o integrals.

 Recall that if $h\in L^{2}([0,1] )$ with $\Vert h\Vert _{L^{2}([0,1])}=1$ then (see e.g. \cite{HN})
\begin{equation*}
sign (W(h))= \sum _{k\geq 0}  b_{2k+1}I_{2k+1}(h^{\otimes (2k+1)} ) \mbox{ with } b_{2k+1}= \frac{2(-1)^{k} }{\sqrt{2\pi }(2k+1) k! 2^{k}}, \hskip0.2cm k\geq 0.
\end{equation*}

We have
\begin{prop}\label{pex}
The standard normal random variable $X$ given by (\ref{x}) is not equal a.s. to $\beta _{1}$ where $\beta$ is its associated DDS martingale.
\end{prop}
{\bf Proof: } By the product formula (\ref{product}) we can express $X$ as (note that $h_{1}$ and $h_{2}$ are orthogonal and there are not contractions of order $l\geq 1$)
\begin{eqnarray*}
&& X= \sum_{k\geq 0}b_{2k+1}I_{2k+2} \left( h_{1} \tilde{ \otimes } h_{2} ^{\otimes 2k+1} \right) \mbox{ and }\\
 &&\mathbf{E}\left( X|{\cal{F}}_{t}\right) =\sum_{k\geq 0}b_{2k+1}I_{2k+2} \left( (h_{1} \tilde{ \otimes } h_{2} ^{\otimes 2k+1})1_{[0,t]}^{\otimes 2k+2}(\cdot ) \right) \mbox{ for every } t\in [0,1].
\end{eqnarray*}
We have
\begin{equation}\label{hat}
(h_{1}\tilde{\otimes } h_{2} ^{\otimes 2k+1})(t_{1},\ldots , t_{2k+2}) = \frac{1}{2k+2} \sum _{i=1} ^{2k+1}h_{1}(t_{i})h_{2}^{\otimes 2k+1}(t_{1}, .., \hat{t}_{i}, .., t_{2k+2})
\end{equation}
where $\hat{t}_{i}$ means that the variable $t_{i}$ is missing. Now, $M^{X}_{t}= \mathbf{E}\left( X|{\cal{F}}_{t}\right)=\int_{0}^{t}u_{s}dW_{s} $ where, by (\ref{hat})
\begin{eqnarray*}
u_{s}&=& \sum_{k\geq 0}b_{2k+1} (2k+2) I_{2k+1}\left( (h_{1} \tilde{\otimes } h_{2}^{2k+1}) (\cdot ,s ) 1_{[0,s]}^{\otimes 2k+1 } (\cdot )\right)\\
&=& \sum_{k\geq 0}b_{2k+1}\left[  h_{1}(s) I_{2k+1}\left( h_{2}^{\otimes 2k+1} 1_{[0,s]} ^{\otimes 2k+1}(\cdot) \right) \right.   \\
&&\left. + (2k+1) h_{2}(s) I_{1}(h_{1}1_{[0,s]}(\cdot ))I_{2k}\left( h_{2} ^{\otimes 2k}1_{[0,s]}^{\otimes 2k}(\cdot ) \right).\right]
\end{eqnarray*}
for every $s\in [0,1]$. Note first that, due to the choice of the functions $h_{1}$  and $h_{2}$,
$$h_{1}(s)  h_{2}(u) 1_{[0,s] }(u)=0 \mbox{ for every } s,u \in [0,1].$$
Thus  the first summand of $u_{s}$ vanishes and
\begin{equation*}
u_{s}= \sum_{k\geq 0}b_{2k+1} (2k+1) h_{2}(s) I_{1}(h_{1}1_{[0,s]}(\cdot ))I_{2k}\left( h_{2} ^{\otimes 2k}1_{[0,s]}^{\otimes 2k}(\cdot ) \right).
\end{equation*}
Note also that $h_{1}(x)1_{[0,s]}(x)=h_{1}(x)$  for every $s$ in the interval $[\frac{1}{2}, 1]$. Consequently, for every $s\in [0,1]$
\begin{equation*}
u_{s}= W(h_{1}) \sum_{k\geq 0}b_{2k+1} (2k+1) h_{2}(s) I_{2k}\left( h_{2} ^{\otimes 2k}1_{[0,s]}^{\otimes 2k}(\cdot ) \right).
\end{equation*}

Let us compute the chaos decomposition of the random variable $\int_{0}^{1} u_{s}^{2}ds.$
Taking into account the fact that $h_{1}$ and $h_{2}$ have disjoint support we can write
\begin{eqnarray*}
&&\int_{0}^{1} u_{s}^{2}ds\\
&=&  \sum_{k,l\geq 0}b_{2k+1}b_{2l+1}(2k+1)(2l+1)W(h_{1})^{2} \int_{0}^{1}ds h_{2}(s)^{2} I_{2k}\left( h_{2} ^{\otimes 2k}1_{[0,s]}^{\otimes 2k}(\cdot ) \right)I_{2l}\left( h_{2} ^{\otimes 2l}1_{[0,s]}^{\otimes 2l}(\cdot ) \right).
\end{eqnarray*}
Since
\begin{equation*}
W(h_{1})^{2} = I_{2} \left( h_{1}^{\otimes 2} \right) + \int_{0}^{1} h_{1}(u)^{2}du=I_{2} \left( h_{1}^{\otimes 2} \right)+1
 \end{equation*}
and
\begin{equation*}
\mathbf{E} \left( sign (W(h_{2}) \right) ^{2}=\int_{0}^{1}ds h_{2}^{2}(s)\mathbf{E} \left( \sum_{k\geq 0}b_{2k+1}(2k+1) I_{2k}  \left( h_{2} ^{\otimes 2k}1_{[0,s]}^{\otimes 2k}(\cdot )\right) \right) ^{2}=1
\end{equation*}
we get
\begin{eqnarray*}
&&\int_{0}^{1} u_{s}^{2}ds=\left( 1+ I_{2} \left( h_{1}^{\otimes 2} \right)\right) \\
&&\times \left( 1+ \sum_{k,l\geq 0}b_{2k+1}b_{2l+1}(2k+1)(2l+1) \int_{0}^{1}ds h_{2}(s)^{2}\right. \\
&& \left. \left[ I_{2k}\left( h_{2} ^{\otimes 2k}1_{[0,s]}^{\otimes 2k}(\cdot ) \right)I_{2l}\left( h_{2} ^{\otimes 2l}1_{[0,s]}^{\otimes 2l}(\cdot ) \right)-\mathbf{E} I_{2k}\left( h_{2} ^{\otimes 2k}1_{[0,s]}^{\otimes 2k}(\cdot ) \right)I_{2l}\left( h_{2} ^{\otimes 2l}1_{[0,s]}^{\otimes 2l}(\cdot ) \right)\right] \right)\\
&=:& \left( 1+ I_{2} \left( h_{1}^{\otimes 2} \right)\right)(1+A).
\end{eqnarray*}
Therefore we obtain that $\int_{0}^{1}u_{s}^{2}ds =1$ almost surely if and only if $\left( 1+ I_{2} \left( h_{1}^{\otimes 2} \right)\right)(1+A)=1$ almost surely which implies that  $I_{2}(h_{1}^{\otimes 2} )(1+A) +A=0$ a.s. and this is impossible because $I_{2}(h_{1}^{\otimes 2})$ and $A$ are independent. \qed

\vskip0.5cm

We obtain an interesting consequence of the above result.

\begin{corollary}
Let $X$ be given by (\ref{x}). Then the bracket of the martingale
$M^{X}$ with $M^{X}_{t}=\mathbf{E}\left(X | {\cal{F}}_{t}\right)$ is
not bounded by 1.
\end{corollary}
{\bf Proof: } It is a consequence of Proposition \ref{pex} and of Theorem \ref{tt}.\qed

\vskip0.3cm

\begin{remark}
Proposition \ref{pex} provides an interesting example of a Brownian motion $\beta$ and of a stopping time $T$ for its filtration such that
$\beta _{T}$ is standard normal and $T$ is not almost surely equal to 1. \end{remark}

 Let us make a short summary  of the results in the first part of our paper: if $X$ is a standard normal random variable and
 the bracket of $M^{X}$ is bounded a.s. by 1 then $X$ can be expressed almost surely
  as a Wiener integral with respect to a Brownian motion  on the same (or possibly extended)
   probability space. The Brownian is obtained via DDS theorem. The property is still true when the bracket is
    bounded and $T$ and $\beta _{T}$ are independent random variables. If the bracket of $M^{X}$ is not bounded,
    then $X$ is not necessarily  equal with $\beta _{1}$, $\beta$ being its associated DDS Brownian motion. This is the case of the variable (\ref{x}).

Nevertheless, we will see that after a suitable extension of the probability space, any standard normal random variable can be written as the value at time 1 of a Brownian motion constructed on this extended probability space.

\begin{prop}\label{peccati}
Let $X_{1}$ be a standard normal random variable on $(\Omega _{1}, {\cal{F}}_{1}, P_{1})$ and for every $i\geq 2$ let $(\Omega _{i}, {\cal{F}}_{i}, P_{i}, X_{i})$ be independent copies of $(\Omega _{1}, {\cal{F}}_{1}, P_{1}, X_{1})$. Let $(\Omega _{0}, {\cal{F}}_{0}, P_{0}) $ be the product probability space. On $\Omega _{0}$ define for every $t\in[0,1]$
\begin{equation*}
W^{0}_{t}= \sum_{k\geq 1}f_{k}(t) X_{k}
\end{equation*}
where $(f_{k})_{k\geq 1}$ are orthonormal elements of
$L^{2}([0,1])$. Then $W^{0}$ is a Brownian motion on $\Omega _{0}$
and $X_{1}=\int_{0}^{1}\left( \int_{u}^{1} dsf_{1}(s)
\right)dW^{0}_{u}$ a.s..
\end{prop}
{\bf Proof: } The fact that $W^{0}$ is a Brownian motion is  a consequence of the Karhunen-Lo\`eve  theorem. Also, note that
$$X_{1}=\langle W^{0}, f_{1}\rangle = \int_{0}^{1} W^{0}_{s}f_{1}(s)ds$$
and the conclusion is obtained by interchanging the order of integration.  \qed

\vskip0.5cm

\begin{remark}
Let us denote by ${\cal{F}}^{0}_{t}$ the natural filtration of $W^{0}$. It also holds that
$$E\left( X_{1}| {\cal{F}}^{0}_{t} \right)= E \int_{0} ^{t} g_{u} dW^{0}_{u}$$ where $g_{u}= \int_{u}^{1} dsf_{1}(s) $.  It is obvious that
the martingale $E\left( X_{1}| {\cal{F}}^{0}_{t} \right)$ is a Brownian motion via the  DDS  theorem and $X_{1}$ can be expressed as a Brownian at time 1.
\end{remark}

\section{Consequences}
We think that the consequences of this result are multiple. We will prove here first that a random variable $X$  which
 lives in a finite sum of Wiener chaoses cannot be Gaussian if the bracket of $M^{X}$ is bounded by 1.
 Again we fix a Wiener process $(W_{t})_{t\in [0, 1]}$ on $\Omega$.

 Let us start with the following lemma.

\begin{lemma}\label{fc}
Fix $N\ge 1$. Let $g \in L^{2}([0,1] ^{\otimes N+1}) $ symmetric in its first $N$ variables  such that  $\int_{0}^{1} ds g(\cdot ,s ) \tilde{\otimes } g(\cdot , s)=0$ almost everywhere  on $[0,1]^{\otimes 2N}$. Then for every $k=1, \ldots, N-1$ it holds that
\begin{equation*}
\int_{0}^{1} dsg(\cdot , s)\tilde{\otimes }_{k} g(\cdot ,s)= 0 \mbox { a.e. on } [0,1] ^{2N-2k}.
\end{equation*}
\end{lemma}
{\bf Proof: } Without  loss of generality we can assume that $g$ vanish on the diagonals $(t_{i}=t_{j})$ of $[0,1] ^{\otimes (N+1)}$. This is possible from the construction of multiple stochastic  integrals. From the hypothesis,   the function
\begin{equation*}
(t_{1}, \ldots ,t_{2N})\to \frac{1}{(2N)!}\sum_{\sigma \in S_{2N}}\int_{0}^{1}ds g(t_{\sigma (1)}, \ldots , t_{\sigma (N)},s)g(t_{\sigma (N+1)}, \ldots , t_{\sigma (2N)}, s)
\end{equation*}
vanishes almost everywhere on $[0,1]^{\otimes 2N}$. Put $t_{2N-1} =t_{2N}=x\in [0,1]$. Then for every $x$, the function
\begin{equation*}
(t_{1}, \ldots  t_{2N-2} ) \to \sum_{\sigma \in S_{2N-2}}\int_{0}^{1}ds g (t_{\sigma (1)}, \ldots , t_{\sigma (N-1)},x, s)g(t_{\sigma (N)}, \ldots , t_{\sigma (2N-2)},x, s)
\end{equation*}
is zero a.e. on $[0,1]^{\otimes (2N-2)}$ and integrating with respect to $x$ we obtain that $\int_{0}^{1}ds g(\cdot , s)\tilde{\otimes }_{1} g(\cdot , s) =0 $ a.e. on $[0,1]^{\otimes (2N-2)}$. By repeating the procedure we obtain the conclusion.  \qed

\vskip0.3cm
Let us also recall the following result from \cite{NO}.
\begin{prop}\label{pno}
Suppose that $F=I_{N}(f_{N})$ with $f \in L^{2}([0,1] ^{N})$ symmetric  and $N\geq 2$ fixed. Then the distribution of $F$ cannot be normal.
\end{prop}

We are going to prove the same property for variables that can be expanded  into a finite sum of multiple integrals.
\begin{theorem}
Fix $N\geq 1$ and et let $X$ be a centered random variable such that
$X=\sum_{n=1} ^{N+1} I_{n} (f_{n})$ where $f\in L^{2}([0,1] ^{n}) $
are symmetric functions. Suppose that the bracket of the martingale
$M^{X}$ (\ref{mx})  is bounded almost surely by 1. Then the law of
$X$ cannot be normal.
\end{theorem}
{\bf Proof: } We will assume that $\mathbf{E}X^{2} =1$. Suppose that $X$ is standard normal. We can write $X$ as $X=\int_{0}^{1} u_{s}dW_{s} $ where $u_{s} =\sum_{n=1}^{N} I_{n}(g_{n}(\cdot ,s) )$.  As a consequence of Proposition \ref{peccati},
\begin{equation*}
\int_{0}^{1} u_{s}^{2}ds =1 \hskip0.5cm \mbox { a. s. }
\end{equation*}
But from the product formula (\ref{product})
\begin{eqnarray*}
\int_{0}^{1}u_{s}^{2}ds&=&\int_{0}^{1} ds \left(  \sum_{n=1}^{N} I_{n}(g_{n}(\cdot ,s) )\right) ^{2} \\
&=& \int_{0}^{1}ds \sum_{m,n=1}^{N} \sum _{k=1}^{m\wedge n} k! C_{n}^{k} C_{m}^{k} I_{m+n-2k} (g_{n}(\cdot , s)\otimes g_{m}(\cdot ,s ) ) ds.
\end{eqnarray*}
The idea is to benefit from the fact that the highest order chaos, which appears only once in the above  expression, vanishes. Let us look to the chaos of  order $2N$ in the above decomposition. As we said,  it appears only  when we multiply $I_{N} $ by $I_{N}$ and consists in the random variable $I_{2N} \left( \int_{0}^{1} g_{N} (\cdot ,s ) \otimes g_{N} (\cdot ,s) ds\right)$. The isometry of multiple integrals (\ref{iso}) implies that
 \begin{equation*}
 \int_{0}^{1} g_{N} (\cdot ,s ) \tilde{\otimes } g_{N} (\cdot ,s) ds =0 \mbox { a. e. on } [0,1] ^{2N}
 \end{equation*}
and by Lemma \ref{fc}, for every $k=1, \ldots , N-1$.
\begin{equation}\label{t100}
 \int_{0}^{1} g_{N} (\cdot ,s ) \tilde{\otimes }_{k} g_{N} (\cdot ,s) ds =0 \mbox { a. e. on } [0,1] ^{2N-2k}.
\end{equation}
Consider now the the random variable $Y:= I_{N+1}(f_{N+1})$.
 It can be written as $Y=\int_{0}^{1} I_{N} (g_{N}(\cdot ,s))dW_{s}$ and  b
y the DDS theorem, $Y=\beta ^{Y}_{\int_{0}^{1} ds (I_{N}(g_{N}(\cdot
,s ) ))^{2}}$.
 The multiplication formula together with (\ref{t100}) shows that $\int_{0}^{1} ds (I_{N}(g_{N}(\cdot ,s ) ))^{2}$ is deterministic and as a consequence $Y$ is Gaussian. This is in contradiction with Proposition \ref{pno}. \qed

\vskip0.2cm
The conclusion of the above theorem still holds if
$M^{X}$ satisfies  (\ref{H}) and $\langle M^{X}\rangle _{1}$ is
independent by $\beta _{ \langle M^{X}\rangle _{1}}$.

\vskip0.2cm

Finally let us make a connection with several recent results obtained via Stein's method and Malliavin calculus. Recall that the Ornstein-Uhlenbeck operator is defined as
$LF=-\sum_{n\geq 0} nI_{n}(f_{n})$ if $F$ is given by (\ref{sum1}). There exists a connection between $\delta, D $ and $L$ in the sense that a random variable $F$ belongs to the domain of $L$ if and only if $F\in \mathbb{D}^{1,2}$ and $DF \in Dom (\delta)$ and then $\delta DF=-LF$.

Let us denote by $D$ the Malliavin derivative with respect to $W$ and let, for any $X\in \mathbb{D}^{1,2}$
\begin{equation*}
G_{X}= \langle DX, D(-L)^{-1} X \rangle .
\end{equation*}

The following theorem is a collection of results in several recent papers.
\begin{theorem}\label{t1}
Let $X$ be a random variable in the space $\mathbb{D}^{1,2}$. Then the following affirmations are equivalent.
\begin{description}
\item{1. } $X$ is a standard normal random variable.
\item{2. } For every $t \in \mathbb{R}$, one has $\mathbf{E}\left( e^{itX} (1-G_{X} ) \right)=0$.
\item{3. } $\mathbf{E}\left( (1-G_{X})/X\right) =0$.
\item{4.} For every $z\in \mathbb{R} $, $\mathbf{E}\left( f'_{z} (1-G_{X}) \right)= 0$, where $ f_{z} $ is the solution of the Stein's equation (see \cite{NoPe1}).
\end{description}
\end{theorem}
{\bf Proof: }We will show that $1. \Rightarrow  2. \Rightarrow 3.\Rightarrow 4. \Rightarrow 1$. First suppose that $X \sim N(0,1)$.  Then
\begin{eqnarray*}
\mathbf{E}\left( e^{itX } (1-G_{X})\right)&=& \mathbf{E}(e^{itX}) -\frac{1}{it} \mathbf{E} \langle De^{itX}, D(-L) ^{-1} X \rangle \\
&=&  \mathbf{E}(e^{itX_{n}}) -\frac{1}{it} \mathbf{E} \left( X e^{itX}\right)=\varphi_{X}(t)-\frac{1}{t} \varphi '_{X} (t)=0.
\end{eqnarray*}
Let us prove now the implication $2. \Rightarrow 3. $ It has also proven in \cite{NV}, Corollary 3.4. Set $F=1-G_{X}$. The random variable $\mathbf{E}(F|X)$ is the Radon-Nykodim derivative with respect to $P$ of the measure $Q(A)= \mathbf{E}(F1_{A})$, $A\in \sigma(X)$. Relation 1. means that
$\mathbf{E}\left( e^{itX} \mathbf{E}(F / X)\right)=\mathbf{E}_{Q } (e^{itX}) =0$  and  consequently $Q (A)=\mathbf{E}(F1_{A}) = 0$  for any $A \in \sigma (X_{n})$. In other words, $\mathbf{E}(F|X) =0.$
The implication $3. \Rightarrow 4$ is trivial and the implication $4. \Rightarrow1. $ is a consequence of a result in \cite{NoPe1}. \qed

\vskip0.5cm

As we said, this property can be easily understood and checked if
$X$ is in the first Wiener chaos with respect to $W$. Indeed, if
$X=W(f)$ with $\Vert f\Vert _{L^{2}([0,1])} =1$ then
$DX=D(-L)^{-1}X=f$ and clearly $G_{X}=1$. There is no need to
compute the conditional expectation given $X$, which is in practice
very difficult to be computed. Let us consider now the case of the
random variable $Y=\int_{0}^{1} sign (W_{s}) dW_{s}$. The chaos
expansion of this variable is known. But $Y$ is not even
differentiable in the Malliavin sense so it is not possible to check
the conditions from Theorem \ref{t1}. Another example is related to
the   Bessel process (see the random variable \ref{bes}). Here again
the chaos expansion of $X$ can be obtained (see e.g. \cite{HN}) but
is it impossible to compute the conditional expectation given $X$.

But on the other hand, for both variables treated above their is another explanation of their normality which comes from L\'evy's characterization theorem.  Another explanation can be obtained from the results in Section 2. Note that  these two examples are random variables such that  the bracket of $M^{X}$ is bounded a.s.

\begin{corollary}
Let $X$ be an integrable random variable on  $(\Omega , {\cal{F}}, P)$. Then $X$ is a standard normal random variable if and only if there exists a Brownian motion $(\beta _{t})_{t\geq 0} $ on an extension of $\Omega$ such that
\begin{equation}\label{bb}
\langle D^{\beta } X, D^{\beta }(-L ^{\beta })^{-1} X\rangle =1.
\end{equation}
\end{corollary}
{\bf Proof: } Assume that $X\sim N(0,1)$. Then by  Proposition \ref{peccati}, $X=\beta _{1}$ where $\beta $ is a Brownian motion on an extended probability space. Clearly (\ref{bb}) holds.
Suppose that there exists $\beta$ a Brownian motion on $(\Omega , {\cal{F}}, P)$ such that (\ref{bb}) holds. Then for any  continuous  and piecewise differentiable function $f$ with $\mathbf{E}f'(Z) <\infty$ we have
\begin{eqnarray*}
\mathbf{E}\left( f'(Z)-f(X)X \right) &=&\mathbf{E}\left( f'(X)-f'(X)\langle D^{\beta } X, D^{\beta }(-L ^{\beta })^{-1} X\rangle\right) \\
&=& \mathbf{E}\left( f'(Z) (1-\langle D^{\beta } X, D^{\beta }(-L ^{\beta })^{-1} X\rangle\right)=0
\end{eqnarray*}
and this implies that $X\sim N(0,1)$ (see \cite{NoPe1}, Lemma 1.2). \qed

\vskip0.5cm

{\bf Acknowledgement: } The Proposition \ref{peccati} has been
introduced in the paper after a discussion with Professor Giovanni
Peccati. We would to thank him for this. We are also grateful to
Professor P. J. Fitzsimmons for detecting a mistake in the first
version of this work.

\end{document}